\documentclass[number,citesort,seceqn,dvips]{arxbj}

\aid{0}
\volume{17}
\issue{2}
\pubyear{2011}
\firstpage{671}
\lastpage{686}
\doi{10.3150/10-BEJ287}

\makeatletter
\newcommand{\R}{\mathbb{R}}
\newcommand{\iint}{\int\!\!\!\int}
\newtheorem{theorem}{Theorem}[section]
\newtheorem{lemma}{Lemma}[section]
\newtheorem{corollary}{Corollary}[section]

\newremark{remark}{Remark}[section]
\newremark{example}{Example}
\makeatother
\begin{document}
\begin{frontmatter}

\title{Limit theorems for functions of~marginal~quantiles}
\runtitle{Functions of marginal quantiles}

\begin{aug}
\author[a]{\fnms{G.~Jogesh} \snm{Babu}\thanksref{a}\ead[label=e1]{babu@stat.psu.edu}},
\author[b]{\fnms{Zhidong} \snm{Bai}\thanksref{b,e2}\ead[label=e2,mark]{stabaizd@nus.edu.sg}},
\author[b]{\fnms{Kwok Pui} \snm{Choi}\corref{}\thanksref{b,e3}\ead[label=e3,mark]{stackp@nus.edu.sg}}\break
\and
\author[c]{\fnms{Vasudevan} \snm{Mangalam}\thanksref{c}\ead[label=e4]{mangalam@fos.ubd.edu.bn}}

\runauthor{Babu, Bai, Choi and Mangalam}

\address[a]{Department of Statistics,
326 Joab L.~Thomas Building,
The Pennsylvania State University,
University Park,
PA 16802-2111,
USA.
\printead{e1}}

\address[b]{Department of Statistics and Applied Probability,
National University of Singapore,
6 Science Drive 2,
Singapore 117546.
\printead{e2,e3}}

\address[c]{Department of Mathematics,
Universiti Brunei Darussalam,
Brunei.\\
\printead{e4}}
\end{aug}

\received{\smonth{12} \syear{2009}}
\revised{\smonth{5} \syear{2010}}

%
\begin{abstract}
Multivariate distributions are explored using the joint distributions
of marginal sample quantiles.
Limit theory for the
mean of a function of order statistics is presented. The results
include a multivariate central limit theorem and a strong law of large
numbers. A result similar to Bahadur's representation of quantiles
is established for the mean of a function of the marginal quantiles.
In particular, it is shown that
\[
\sqrt n \Biggl( \frac{1}{n} \sum_{i=1}^n \phi\bigl(X_{n\dvtx  i }^{(1)},
\ldots,
X_{n\dvtx  i }^{(d)}\bigr) - \bar{\gamma} \Biggr) = \frac{1}{\sqrt{n}} \sum
_{i=1}^nZ_{n, i} +
\mathrm{o}_{P}(1)
\]
as $n \rightarrow\infty$, where $\bar{\gamma}$ is a constant
and $Z_{n, i}$ are i.i.d.~random variables for each $n$.
This leads to the central limit theorem.
Weak convergence to a Gaussian process using equicontinuity of
functions is indicated.
The results are established under very general conditions.
These conditions are shown to be satisfied in many commonly occurring
situations.
\end{abstract}

%
\begin{keyword}
\kwd{central limit theorem}
\kwd{Cram\'{e}r--Wold device}
\kwd{lost association}
\kwd{quantiles}
\kwd{strong law of large numbers}
\kwd{weak convergence of a process}
\end{keyword}

\end{frontmatter}

\section{Introduction}\label{sec1}
Let $\{(X_i^{(1)},X_i^{(2)},\ldots, X_i^{(d)}), i=1,2,\ldots\}$ be
a sequence of random vectors such that for each~$j$ $(1\le j\le
d)$, $\{X_1^{(j)},X_2^{(j)},\ldots\}$ forms a sequence of
independent and identically distributed (i.i.d.)~random
variables. For $1\le j, k \le d$, let $F_j$ and $F_{j,k}$ denote the
distributions of
$X_1^{(j)}$ and $(X_1^{(j)}, X_1^{(k)})$, respectively.
Let $X_{n\dvtx i}^{(j)}$
denote the $i$th order statistic ($\frac{i}{n}$th quantile) of $\{
X_1^{(j)},X_2^{(j)},\ldots,
X_n^{(j)}\}$. The vector $(X_{n\dvtx i}^{(1)},
\ldots, X_{n\dvtx i}^{(d)})$ corresponds to the $i$th marginal order statistics.
In this article, we study the asymptotic behavior of
the mean of a function of marginal sample quantiles:
%
\begin{equation}
\label{intro} \frac{1}{n} \sum_{i=1}^n \phi\bigl(X_{n\dvtx i}^{(1)},
\ldots, X_{n\dvtx i}^{(d)} \bigr)
\end{equation}
as $n \rightarrow\infty$, where $\phi\dvtx  \R^d \rightarrow\R$
satisfies some mild
conditions.

Our results, Theorems \ref{asconvergence} and \ref{asymnormality}
stated below, were motivated in part by one of the authors
considering \cite{VM} the problem of estimating the parameters in a linear
regression model, $Y = \alpha+ \beta X + \epsilon$, when the
linkage between the variables $X$ and $Y$ was either partially or
completely lost. Were the linkage not lost, then the
least-squares estimator for $\beta$ would be given by $
(\sum_{i=1}^n X_i Y_i - n \bar{X_n} \bar{Y_n}
)/\sum_{i=1}^n (X_i-\bar{X_n})^2$, where $\bar{X_n}$ and
$\bar{Y_n}$ denote the sample means of ($X_1, \ldots, X_n$) and
($Y_1, \ldots, Y_n$). When the linkage is lost, a natural candidate
to estimate $\beta$ is the average of this expression over all
possible permutations of the $Y_i$'s. As the term in the denominator
and the second term in the numerator are permutation invariant, it
remains to consider $\frac{1}{n!} \sum_{\pi\in\mathcal{S}_n}
\frac{1}{n} \sum_{i=1}^n X_i Y_{\pi(i)}.$ This expression is bounded
above by $\frac{1}{n} \sum_{i=1}^n X_{n\dvtx  i} Y_{n\dvtx  i}$ and below by
$\frac{1}{n} \sum_{i=1}^n X_{n\dvtx  i} Y_{n\dvtx  n-i+1}$, by the well-known
rearrangement inequality of Hardy--Littlewood--P\'{o}lya (see \cite
{HLP}, Chapter 10). The asymptotic behavior of the lower bound can
be deduced from that of the upper bound. The upper bound,
$\frac{1}{n} \sum_{i=1}^n X_{n\dvtx  i} Y_{n\dvtx i}$, is a special case of
(\ref{intro}). The problem of the loss of association among paired
data has attracted a lot of attention in various contexts, such as the
broken sample problem, file linkage problem and record linkage (see,
e.g., \cite{BH,CH,DG}). See item (3) in
Section \ref{sec4} for further results and a very brief review of the
literature.

We shall first introduce some notation. We shall reserve $\{U_i\}$ for
a sequence of independent random variables distributed
uniformly on $(0,1)$. Let $U_{n\dvtx  i}$ be the $i$th order statistic
of of ($U_1, \ldots, U_n$). For a probability distribution function $F$
and $0 < t <1$, define $F^{-1}(t) = \inf\{x\dvtx  F(x) \ge t\}$.

Let $\phi$ be a real-valued measurable function on $\R^d$. For
$0 <x, x_1, \ldots, x_d < 1$, $\mathbf{x}=(x_1, \ldots, x_d)$, and
$1\le j, k \le d$,
define
\begin{eqnarray}
\psi(\mathbf{x}) & := &\phi(F_1^{-1}(x_1), \ldots
,F_d^{-1}(x_d)), \label{psi} \\
\gamma(x) & :=& \psi(x, x, \ldots, x), \label{gamma} \\
\psi_j(x) & :=& \frac{\partial\psi(\mathbf{x})}{\partial x_j}\bigg|_{(x,
\ldots, x)}, \label{derivatives} \\
\psi_{j, k} (\mathbf{x}) & :=& \frac{\partial^2 \psi(\mathbf
{x})}{\partial x_j\, \partial x_k},\\
\tilde{\psi}_{j, k}(x) & :=& \psi_{j, k} (x, \ldots, x).
\end{eqnarray}

We shall now introduce conditions on $\phi$ that are used in the
results:
\begin{enumerate}[(C1)]
\item[(C1)] The function $\psi(u_1, \ldots, u_d)$ is continuous
at $ u_1 = \cdots= u_d=u,   0<u < 1$. That is, $\psi$ is continuous
at each point on the diagonal of $(0, 1)^d$. The function $\psi$ need
not be bounded.

\item[(C2)] There exist $K$ and $c_0>0$ such that
\[
|\psi(x_1, \ldots, x_d)| \le K\Biggl(1+ \sum_{j=1}^d |\gamma(x_j)|\Biggr)
\
\qquad \mbox{for }  (x_1,\ldots,x_d) \in(0,c_0)^d\cup(1-c_0,1)^d.
\]

\item[(C3)] Let $\mu_{n:i} = i/(n+1)$. For $1\le j, k \le d$,
\[
\frac{1}{n} \sum_{i=1}^{n} \bigl(\mu_{n\dvtx i} ( 1- \mu_{n\dvtx i})
\bigr)^{3/2} (\psi_j (\mu_{n\dvtx i}))^2
\longrightarrow\int_0^1 \bigl(x(1-x)\bigr)^{3/2} (\psi_j(x)
)^2 \,\mathrm{d}x < \infty
\]
and
\[
\frac{1}{n} \sum_{i=1}^{n} \bigl(\mu_{n\dvtx i} ( 1- \mu_{n\dvtx i})
\bigr)^{3/2} |\tilde{\psi}_{j, k} (\mu_{n\dvtx i})|
\longrightarrow\int_0^1\bigl(x(1-x)\bigr)^{3/2} |\tilde{\psi}_{j,
k}(x)| \,\mathrm{d}x < \infty.
\]

\item[(C4)] For all large $m$, there exist $K=K(m) \ge1$ and $\delta>0$
such that
\[
|\psi(\mathbf{y}) -\psi(\mathbf{x}) - \langle\mathbf{y}-\mathbf{x},
\nabla\psi(\mathbf{x}) \rangle|
\le K \sum_{j, k =1}^d |(y_j -x)(y_k -x)| \bigl(1 + |\psi_{j,k} (\mathbf
{x})| \bigr),
\]
whenever
$\Vert\mathbf{y}-\mathbf{x}\Vert_{\ell_1} < \delta$ and
$ {\min_{1 \le j \le d}   y_j(1-y_j) > x(1-x)/m}, $
where   $\mathbf{x}=(x, \ldots, x)$, $\mathbf{y} =(y_1, \ldots, y_d)
\in
(0, 1)^d$. Here,
$\|\mathbf{y}\|_{\ell_1} := |y_1| + \cdots+ |y_d|$ denotes the
$\ell_1$-norm of $\mathbf{y}$ and $\nabla\psi(\mathbf{x})$ denotes the
gradient of $\psi$.
\end{enumerate}
%

Condition (C3) holds if the functions
$(x(1-x))^{3/2} (\psi_j(x))^2$ and
$(x(1-x))^{3/2} |\tilde{\psi} _{j, k}(x)|$
are Riemann integrable over $(0,1)$
and satisfy $K$-pseudo convexity for $1 \le j, k \le d$.
A function $g$ is said to be $K$-pseudo convex if
$g(\lambda x +(1-\lambda) y) \le K (\lambda g(x) + (1-\lambda)g(y))$.

To state the main results, recall the definition of $\gamma$ in (\ref{gamma}).

\begin{theorem}
\label{asconvergence} Let $\{(X_i^{(1)}, X_i^{(2)}, \ldots,
X_i^{(d)}), i=1,2,\ldots\}$ be a sequence of random vectors such
that for each~$j$ $(1\leq j\leq d)$, $\{X_1^{(j)},X_2^{(j)},\ldots\}$
forms a sequence of i.i.d.~random variables. Suppose~$\phi$ satisfies
conditions \textup{(C1)--(C2)}, $F_j$ is continuous for $1\le j \le d$ and
$\gamma$ is Riemann integrable. Then,
\[
\frac{1}{n}\sum_{i=1}^n\phi\bigl(X_{n\dvtx i}^{(1)},\ldots,X_{n\dvtx i}^{(d)}\bigr)
\mathop{\longrightarrow}^{a.s.} \bar{\gamma}
\]
as $n \rightarrow\infty$, where $\bar{\gamma} =
\int_0^1 \gamma(y)\,\mathrm{d}y$.
\end{theorem}

Note that we need only the independence of the $j$th marginal random variables,
for each $j$. The result does not depend on the joint distribution of
$(X_1^{(1)}, \ldots, X_1^{(d)}).$

\begin{theorem}
\label{asymnormality} Let $\mathbf{X}_i = (X_i^{(1)}, \ldots,
X_i^{(d)})$ be i.i.d.~random vectors. Suppose $\phi$ satisfies
conditions \textup{(C1)--(C4)}, $F_j$ is continuous for $1\le j \le d$ and
$\gamma$ is Riemann integrable. Then,
%
\begin{equation}
\label{Bahadur}
\frac{1}{\sqrt{n}} \sum_{i=1}^n \phi\bigl(X_{n\dvtx  i }^{(1)}, \ldots,
X_{n\dvtx  i }^{(d)}\bigr) - \sqrt{n} \bar{\gamma} =
\frac{1}{\sqrt{n}} \sum_{\ell=1}^n Z_{n, \ell} + \mathrm{o}_P(1),
\end{equation}
where
$Z_{n, \ell} = \frac{1}{n} \sum_{i=1}^n \sum_{j=1}^d W_{j, \ell} (i/n)
\psi_j (i/(n+1))$,
$ W_{j, \ell} (x) =I(U_{\ell}^{(j)} \le x) -x$
for $ 1 \le\ell\le n$ and $\bar{\gamma}$ is defined as in Theorem
\ref{asconvergence}.
Further, as $n \rightarrow\infty$,
%
\begin{equation}
\label{eq1}
\frac{1}{\sqrt{n}} \sum_{i=1}^n \phi\bigl(X_{n\dvtx  i }^{(1)}, \ldots,
X_{n\dvtx  i }^{(d)}\bigr) - \sqrt{n} \bar{\gamma} \mathop{\longrightarrow}^{dist}
N(0, \sigma^2),
\end{equation}
where $G_{j, k}(x, y) = F_{j, k}(F_j^{-1}(x),
F_k^{-1}(y))$ and
\begin{eqnarray*}
\sigma^2 = \lim_{n \rightarrow\infty} \operatorname{Var}(Z_{n, 1})
&=& 2 \sum_{j=1}^d\int_0^1 \!\! \int_0^y x(1-y) \psi_j(x) \psi_j(y) \,\mathrm{d}x \,\mathrm{d}y \\
&&{} + 2 \sum_{1 \le j < k \le d} \int_0^1\!\! \int_0^1 \bigl(G_{j,k}(x,
y)-xy\bigr) \psi_j(x) \psi_k(y) \,\mathrm{d}x \,\mathrm{d}y.
\end{eqnarray*}
\end{theorem}

This theorem can be extended to $m$ functions $\phi_1, \ldots,
\phi_m$ simultaneously using the Cram\'{e}r--Wold device (see
\cite{Bill}), as in the corollary below. Let $\psi_j(x; r)$ denote
the partial derivative of $\phi_r(F_1^{-1}(x_1), \ldots,
F_d^{-1}(x_d))$ with respect to $x_j$ evaluated at $x_1 =
\cdots=x_d=x$.

\begin{corollary}
Let $\phi_1, \ldots, \phi_{m}$ satisfy conditions \textup{(C1)--(C4)}. For
$1 \le r \le m$, if we define $T_n (\phi_r) = \sum_{i=1}^n \phi_r
(X_{n\dvtx  i }^{(1)}, \ldots, X_{n\dvtx  i }^{(d)})$ and
$\bar{\gamma}_r = E\phi_r(F_1^{-1}(U), F_2^{-1}(U), \ldots,
F_d^{-1}(U)$, then
\[
\frac{1}{\sqrt{n}} (T_n(\phi_1), \ldots, T_n(\phi_m) ) - \sqrt{n} (\bar
{\gamma}_1, \ldots, \bar{\gamma}_m)
\mathop{\longrightarrow}^{dist} N(0, \Sigma)  \qquad \mbox{as } n \rightarrow\infty,
\]
where the $(r, s)$th element $\sigma_{r,s}$ of $\Sigma$, is given by
\begin{eqnarray*}
&&\sum_{j=1}^d  \int_0^1 \!\!\int_0^y x(1-y) \bigl(\psi_{j}(x; r) \psi
_{j}(y; s)
+ \psi_{j}(x; s) \psi_{j}(y; r)\bigr) \,\mathrm{d}x \,\mathrm{d}y \\
&&\quad{} + \sum_{1 \le j < k \le d} \int_0^1\!\! \int_0^1 \bigl(G_{j,k}(x,
y)-xy\bigr) \bigl(\psi_j(x; r) \psi_k (y; s)
+ \psi_j (x; s) \psi_{k}(y; r)\bigr) \,\mathrm{d}x \,\mathrm{d}y.
\end{eqnarray*}
\end{corollary}

\begin{pf} Use the Cram\'{e}r--Wold device and Theorem
\ref{asymnormality}. In computing $\sigma_{r,s}$, we used
\[
2 \sigma_{r,s} = \lim_{n\rightarrow\infty} \bigl(\operatorname{Var}(Z_{n, 1, r} +
Z_{n, 1, s}) -
\operatorname{Var}(Z_{n, 1, r}) - \operatorname{Var}(Z_{n, 1, s})\bigr),
\]
where $Z_{n, 1, r}= \frac{1}{n} \sum_{i=1}^n \sum_{j=1}^d W_{j, 1}
(i/n) \psi_j(i/(n+1); r).$
\end{pf}

Our results can be adapted to provide a suitable test statistic for
testing equality of marginal distributions against various
alternative hypotheses using suitable choices for~$\phi$.

\begin{remark} Since the finite-dimensional distributions
converge to multivariate normal
distributions, the weak convergence to a Gaussian process indexed by $t
\in T$
($T$ being an interval of $\R$) can be established under a condition
such as equicontinuity of
$\{\phi_t\dvtx  t \in T\}$.
\end{remark}

\begin{remark}
In Theorem \ref{asconvergence}, we just require i.i.d.~for
each component. No further assumptions are made on how the components
are related. We need a stronger assumption in Theorem
\ref{asymnormality}, namely, that the rows are i.i.d.~random
vectors. Interestingly, the variance of the limiting normal only
depends on the 2-dimensional marginal distributions.
\end{remark}

\begin{remark} Conditions (C1) and (C2) are, in general,
easy to
verify. Condition (C3) is used to control the behavior of the function
$\psi$ around the
neighborhood of $(0, \ldots, 0)$ and $(1, \ldots, 1)$ in $(0, 1)^d$.
For example, if we suppose that $X_1^{(j)}$ is uniformly distributed
over $(0,1)$ for $ j=1, 2$ and
$\phi(x, y)\dvtx  =((x+y)/2)^{-\alpha} (1- (x+y)/2)^{-\alpha}$, then (C3) holds
if $0 < \alpha< 1/4$. However, the first limit in (C3) fails if $
\alpha\ge1/4$ and
the second limit in (C3) fails if $\alpha\ge1/2$.
\end{remark}

\begin{remark}\label{rem1.4} By a compactness argument, condition (C1) is
shown to
be equivalent~to
\begin{enumerate}[(C1$'$)]
\item[(C1$'$)] For any $c \in(0,\frac{1}{2}), \lim_{\delta\rightarrow
0}\omega(c,\delta)=0$, where
%
\begin{eqnarray}
\label{C1'}
&&\omega(c,\delta):=\sup\{|\psi(x_1, \ldots, x_d) - \gamma(y)|\dvtx  |x_i-y|<
\delta,   c <y,
\nonumber
\\[-8pt]
\\[-8pt]
\nonumber
&&\hspace*{162pt} x_i <1-c,   1\le i\le d\}.
\end{eqnarray}
\end{enumerate}
\end{remark}

Proofs of Theorems
\ref{asconvergence} and \ref{asymnormality} are given in Sections
\ref{sec2}
and \ref{sec3}, respectively. The results are illustrated by means of examples
and counterexamples
in the last section.

\section{\texorpdfstring{Proof of Theorem \protect\ref{asconvergence}}{Proof of Theorem 1.1}}
\label{sec2}
The main idea of the proof of Theorem \ref{asconvergence} comes from
the observation that
\[
\frac{1}{n}\sum_{i=1}^n\phi\bigl(X_{n\dvtx i}^{(1)}, \ldots,
X_{n\dvtx i}^{(d)}\bigr) =\frac{1}{n}\sum_{i=1}^n
\psi\bigl(U_{n\dvtx i}^{(1)}, \ldots,U_{n\dvtx i}^{(d)} \bigr)
\approx\int_0^1 \psi(u, \ldots, u) \,\mathrm{d}u.
\]
The cases where $i$ is close to 1 or $n$ need to be carefully analyzed
as $\psi$ could be unbounded near~$0$ and $1$.

\begin{pf*}{Proof of Theorem \protect\ref{asconvergence}} Let
$U_i^{(j)}=F_j(X_i^{(j)})$ for $1 \le i \le n, 1 \le j
\le d$. Therefore, $\{U_1^{(j)}, U_2^{(j)},\ldots\}$ forms a sequence of
i.i.d.~uniformly distributed random variables and
$F^{-1}_j(U_i^{(j)})=X_i^{(j)}$ with
probability 1. Recall that $U_{n\dvtx i}^{(j)}$ denotes the $i$th
order statistic of $U_1^{(j)}, \ldots, U_n^{(j)}$. We write $\mu_{n\dvtx  i}
=EU_{n\dvtx i}^{(j)} = i/(n+1)$.
Recall, also, that $\psi(x_1, \ldots, x_d)=\phi(F_1^{-1}(x_1),
\ldots, F_d^{-1}(x_d))$ and that $\gamma(x) = \psi(x, \ldots, x)$.
For any $\epsilon\in(0, c_0)$,
%
\begin{eqnarray}
\frac{1}{n}\sum_{i=1}^n\phi\bigl(X_{n\dvtx i}^{(1)}, \ldots,
X_{n\dvtx i}^{(d)}\bigr)
= \frac{1}{n}\sum_{i=1}^n \psi\bigl(U_{n\dvtx i}^{(1)}, \ldots,U_{n\dvtx i}^{(d)}
\bigr)
= \Gamma_n+R_{n, 1} + R_{n, 2} + R_{n, 3}
\label{main}
\end{eqnarray}
almost surely, where
\begin{eqnarray*}
\Gamma_n &=& \frac{1}{n}\sum_{i=1}^n\gamma(\mu_{n\dvtx  i}), \\
R_{n, 1} &=& \frac{1}{n}\sum_{1\leq i< \epsilon
n}\bigl(\psi\bigl(U_{n\dvtx i}^{(1)}, \ldots,U_{n\dvtx i}^{(d)} \bigr)
-\gamma(\mu_{n\dvtx  i})\bigr), \\
R_{n, 2} &=& \frac{1}{n}\sum_{\epsilon n\leq i\leq
(1-\epsilon) n}\bigl (\psi\bigl(U_{n\dvtx i}^{(1)}, \ldots
,U_{n\dvtx i}^{(d)} \bigr)
-\gamma(\mu_{n\dvtx  i}) \bigr), \\
R_{n, 3} &=& \frac{1}{n}\sum_{(1-\epsilon) n< i\leq n}
\bigl(\psi\bigl(U_{n\dvtx i}^{(1)}, \ldots,U_{n\dvtx i}^{(d)} \bigr) -\gamma
(\mu_{n\dvtx  i}) \bigr).
\end{eqnarray*}
Since $\gamma$ is Riemann integrable, the Riemann sum
\[
\Gamma_n \rightarrow
\int_0^1\gamma(y)\,\mathrm{d}y =E(\phi(F_1^{-1}(U), \ldots, F_d^{-1}(U)
))  \qquad \mbox{as } n\rightarrow\infty.
\]
Thus, it remains to show that $R_{n, i} \mathop{\longrightarrow}^{\mathrm{a.s.}}0$
as $n\rightarrow\infty$ for $i=1, 2$ and $3$.

\smallskip

For $1 \le j \le d$, by then Glivenko--Cantelli lemma, $\sup_{x \in(0,
1)} |\hat{F}_{n; j}(x) - x| \mathop{\longrightarrow}^{\mathrm{a.s.}}0$ as
$n\rightarrow\infty$, where $\hat{F}_{n; j}$ is the empirical
distribution function of $\{U_i^{(j)}\dvtx  1 \le i \le n \}$. For
$1\leq i\leq n, 1 \le j \le d$, we have
\[
\bigl|U_{n\dvtx i}^{(j)} - \mu_{n\dvtx  i} \bigr| \le \biggl|U_{n\dvtx i}^{(j)} - \frac{i}{n}
\biggr| + \frac{1}{n} = \bigl|U_{n\dvtx i}^{(j)} - \hat{F}_{n;
j}\bigl(U_{n\dvtx i}^{(j)}\bigr) \bigr| +\frac{1}{n} \leq\frac{1}{n} + \sup_{x\in(0,
1)} |x - \hat{F}_{n; j}(x)|.
\]
Hence, it follows that as $n \rightarrow\infty$,
%
\begin{equation}
\delta_n: =\max\bigl\{ \bigl|U_{n\dvtx i}^{(j)} - \mu_{n\dvtx  i}\bigr|, 1\le i \le n, 1 \le
j \le d \bigr\} \mathop{\longrightarrow}^{\mathrm{a.s.}} 0.
\label{deltan}
\end{equation}
Recall the definition of $\omega(c, \delta)$ in (\ref{C1'}). Since
$U_{n\dvtx i}^{(j)} \in(\mu_{n\dvtx  i} - \delta_n, \mu_{n\dvtx  i} +
\delta_n)$ for $ 1\le j \le d$ and for each integer
$i$ in the interval $[n\epsilon, n(1-\epsilon)]$, we have
$|\psi(U_{n\dvtx i}^{(1)}, \ldots, U_{n\dvtx i}^{(d)})-\gamma(\mu_{n\dvtx  i})|
\le \omega(\epsilon, \delta_n)$, provided $\delta_n < \epsilon/2$.
Hence, if $\delta_n < \epsilon/2$, by (\ref{deltan}) and (C1$'$) (which
is equivalent to (C1)
by Remark~\ref{rem1.4} in Section \ref{sec1}), we have
\[
|R_{n, 2}| \le\frac{1}{n}\sum_{\epsilon n\le i\le
(1-\epsilon) n} \bigl|\psi\bigl(U_{n\dvtx i}^{(1)}, \ldots,
U_{n\dvtx i}^{(d)}\bigr) -\gamma(\mu_{n\dvtx  i}) \bigr| \le
\omega(\epsilon,\delta_n) \mathop{\longrightarrow}^{\mathrm{a.s.}} 0
\]
as $n\rightarrow\infty$. By (C2),
\[
|R_{n, 1}| \le K \sum_{j=1}^d R_{n, 1, j} +
\frac{1}{n} \sum_{1 \le i < \epsilon n} |\gamma(\mu_{n\dvtx  i})| + K\epsilon
,
\]
where
$R_{n, 1, j} = n^{-1} \sum_{1\le i < \epsilon n}|\gamma
(U_{n\dvtx i}^{(j)})|$ for $1 \le j \le d$.
Clearly, if $U_{n\dvtx  (\epsilon n)+1}^{(j)} \le2 \epsilon$, then
\[
R_{n, 1, j} \leq n^{-1} \sum_{1\leq i\leq n} \bigl|\gamma\bigl(U_i^{(j)}\bigr)\bigr| I\bigl(U_i^{(j)} \leq2\epsilon\bigr).
\]
Note that, with probability 1, $U_{n\dvtx  (\epsilon n)+1}^{(j)} \le2
\epsilon$
for all large $n$ and the right-hand side of the above inequality
goes to $ \int_0^{2\epsilon}|\gamma(y)|\,\mathrm{d}y$ a.s.~as $n \rightarrow\infty$.
Hence,
\[
\limsup_{n\rightarrow\infty}|R_{n, 1}| \le(Kd+1) \biggl(\epsilon+ \int
_0^{2\epsilon}|\gamma(y)|\,\mathrm{d}y \biggr) \qquad   \mbox{a.s.}
\]
As $|\gamma|$ is integrable, letting $\epsilon$ tend to zero, we
conclude that $R_{n, 1} \mathop{\longrightarrow}^{\mathrm{a.s.}} 0$. A similar \mbox{argument} will
show that $R_{n, 3} \mathop{\longrightarrow}^{\mathrm{a.s.}} 0$ as $n \rightarrow\infty$. This
completes the proof of Theorem~\ref{asconvergence}.
\end{pf*}

\section{\texorpdfstring{Proof of Theorem \protect\ref{asymnormality}}{Proof of Theorem 1.2}}
\label{sec3}

As in the proof of Theorem \ref{asconvergence}, we
introduce $U_{i}^{(j)} = F_j(X_i^{(j)})$ for $ 1 \le i \le n, 1
\le j \le d$. It follows that $(U_i^{(1)}, \ldots, U_i^{(d)}), 1
\le i \le n$, are i.i.d.~random vectors. For $1 \le j, k \le d$,
note that~$G_{j, k}$ is the joint distribution of $(U_1^{(j)}, U_1^{(k)})$.
In particular, $G_{j, j}(x, y) = \min\{x, y\}$, for $1 \le j \le d$.
Using the notation introduced in Section \ref{sec1}, we outline some key
approximations used in the proof
of Theorem \ref{asymnormality}. In particular, (\ref{Bahadur}) follows from
\begin{eqnarray*}
&&\frac{1}{\sqrt{n}} \sum_{i=1}^n \psi\bigl(U_{n\dvtx  i}^{(1)}, \ldots,
U_{n\dvtx  i}^{(d)}\bigr) - \sqrt{n} \bar{\gamma}\\
&&\quad \approx
\frac{1}{\sqrt{n}} \sum_{i=1}^n \bigl(\psi\bigl(U_{n\dvtx  i}^{(1)}, \ldots,
U_{n\dvtx  i}^{(d)}\bigr) - \gamma(\mu_{n, i})\bigr) \\
&&\quad\approx
\frac{1}{\sqrt{n}} \sum_{i=1}^n \sum_{j=1}^d \bigl( U_{n\dvtx  i}^{(j)} - \mu
_{n,i}\bigr) \psi_j(\mu_{n\dvtx i}) \\
&&\quad\approx \frac{1}{\sqrt{n}} \sum_{i=1}^n \sum_{j=1}^d \sum_{\ell=1}^n
\bigl(I\bigl(U_{\ell}^{(j)} - i/n\bigr)\bigr) \psi_j(\mu_{n\dvtx i}).
\end{eqnarray*}
The proof of the first approximation, which is about $\sqrt{n}$ times
the difference between the Riemann sum and the integral $\bar{\gamma}$, is
non-trivial and is handled in Lemma \ref{lemma4}. We use Bahadur's
representation of quantiles in the last approximation.
We start with some technical lemmas,
the first of which is well known (see \cite{David}, page 36).

\begin{lemma}
\label{lemma1} Suppose that $U_{n\dvtx 1} \le\cdots\le U_{n\dvtx n}$ denote the order
statistics of $n$ independent random variables that are uniformly distributed
over $(0, 1)$. Then, for $ 1\le i \le n$,
\[
\operatorname{Var} (U_{n\dvtx i}) = \frac{\mu_{n\dvtx i} (1-\mu_{n\dvtx  i})}{(n+2)} \le\frac{1}{n}.
\]
\end{lemma}

\begin{lemma}
\label{lemma2} Under condition \textup{(C3)}, the limiting variance $\sigma^2$
is well defined.
\end{lemma}

\begin{pf}It suffices to show that for $ 1 \le j, k \le d$,
\begin{eqnarray}
\beta_1&:=& \iint_{0 < x < y <1} |G_{j, k} (x, y)-xy| |\psi_j(x) \psi
_k (y)| \,\mathrm{d}x \,\mathrm{d}y < \infty, \label{bt1} \\
\beta_2&:=& \iint_{0 < y < x <1} |G_{j, k} (x, y)-xy| |\psi_j(x)
\psi_k (y)| \,\mathrm{d}x \,\mathrm{d}y < \infty. \label{bt2}
\end{eqnarray}
To prove (\ref{bt1}), we introduce $W_j(x) := I(U_1^{(j)} \le x)
-x$. Here, $W_j(x)$ has mean 0 and variance $x(1-x)$. Furthermore,
$EW_j(x) W_k(y) = G_{j,k}(x, y) -xy.$
Thus, $EW_j(x)W_j(y)=x(1-y)$ when $x<y$. By the Cauchy--Schwarz
inequality, $\beta_1^2$ is bounded above by
\begin{eqnarray*}
&&\biggl(E \int_0^1\!\! \int_0^y  \bigl(x/(1-y)\bigr)^{1/4}
|\psi_j(x)| |W_j(x)| \bigl((1-y)/x\bigr)^{1/4} |\psi_k(y)| |W_k(y)| \,\mathrm{d}x \,\mathrm{d}y
\biggr)^2 \\
&&\quad\le E \int_0^1 \!\!\int_0^y \bigl(x/(1-y)\bigr)^{1/2} (\psi_j(x))^2 (W_j(x))^2 \,\mathrm{d}x
\,\mathrm{d}y \\
&&\hspace*{34pt}\qquad{}  \times E \int_0^1 \!\!\int_0^y \bigl((1-y)/x\bigr)^{1/2} (\psi_k(y))^2
(W_k(y))^2 \,\mathrm{d}x \,\mathrm{d}y \\
&&\quad= \int_0^1\!\! \int_0^y x^{3/2} (1-x) (1-y)^{-1/2} (\psi_j(x))^2 \,\mathrm{d}x \,\mathrm{d}y \\
&&\hspace*{25pt}\qquad{}  \times\ \int_0^1 \!\!\int_0^y x^{-1/2}y(1-y)^{3/2} (\psi_k(y))^2\,\mathrm{d}x \,\mathrm{d}y \\
&&\quad= 4 \int_0^1 x^{3/2}(1-x)^{3/2} (\psi_j(x))^2 \,\mathrm{d}x\\
&&\hspace*{16pt}{}\times \int_0^1
y^{3/2}(1-y)^{3/2} (\psi_k(y))^2 \,\mathrm{d}x < \infty.
\end{eqnarray*}
Similarly, we can prove (\ref{bt2}). This completes the proof of
Lemma \ref{lemma2}.
\end{pf}

\begin{lemma}
\label{lemma4} Let $\phi\dvtx  (0, 1)^d \to\R$ satisfy
condition \textup{(C3)}.
Suppose that the function $\gamma$ associated with $\phi$ and defined
in (\ref{gamma}) is Riemann
integrable. We then have
\[
\frac{1}{\sqrt{n}} \sum_{i=1}^{n} \gamma(\mu_{n\dvtx  i})
-\sqrt{n} \int_0^1 \gamma(x) \,\mathrm{d}x \to0
\]
as $ n \rightarrow\infty.$
\end{lemma}

\begin{pf}As $\gamma'(x) = \psi_1(x) + \cdots+ \psi_d(x)$, condition
(C3) implies that $(x(1-x))^{3/2} (\gamma'(x))^2$ is Riemann
integrable. We have
\begin{eqnarray*}
&&\frac{1}{\sqrt{n}} \sum_{i=1}^{n} \gamma(\mu_{n\dvtx  i})
-\sqrt{n} \int_0^1 \gamma(x) \,\mathrm{d}x\\
&&\quad= \sqrt{n} \sum_{i=1}^n \int_{(i-1)/n}^{i/n}
\bigl(\gamma(\mu_{n\dvtx  i}) -\gamma(x)\bigr) \,\mathrm{d}x \\
&&\quad= \sqrt{n} \sum_{i =1}^n \biggl( \int_{(i-1)/n}^{\mu_{n\dvtx  i}}
\!\!\int_x^{\mu_{n\dvtx  i}} \gamma'(y) \,\mathrm{d}y \,\mathrm{d}x
- \int_{\mu_{n\dvtx  i}}^{i/n} \!\!\int_x^{\mu_{n\dvtx  i}} \gamma'(y) \,\mathrm{d}y \,\mathrm{d}x \biggr)\\
&&\quad= \sqrt{n} \int_{0}^{1} g_n(y) \gamma'(y) \,\mathrm{d}y,
\end{eqnarray*}
where
\[
g_n(y) =
\cases{
y-(i-1)/n, & \quad $\mbox{if } (i-1)/n \le y < i/(n+1), 1 \le i \le n,$
\vspace*{2pt}\cr
y- i/n, & \quad $\mbox{if } i/(n+1) \le y < i/n, 1 \le i \le n.$}
\]
Note that
\[
|g_n(y)| \le
\cases{
y, & \quad $\mbox{if } 0 <y \le1/n,$ \vspace*{2pt}\cr
1/n, & \quad $\mbox{if } 1/n < y < 1- 1/n,$ \vspace*{2pt}\cr
1-y, & \quad$\mbox{if } 1 - 1/n \le y <1.$}
\]
Therefore,
\begin{eqnarray*}
&&\Biggl(\frac{1}{\sqrt{n}} \sum_{i=1}^{n} \gamma(\mu_{n\dvtx  i})
-\sqrt{n} \int_0^1 \gamma(x) \,\mathrm{d}x\Biggr)^2 \\
&&\quad= n \biggl( \int_0^1 g_n(y) \gamma'(y) \,\mathrm{d}y\biggr)^2 \\
&&\quad\le n \int_0^1 (g_n(y))^2 \bigl(y(1-y)\bigr)^{-3/2} \,\mathrm{d}y \int_0^1
\bigl(y(1-y)\bigr)^{3/2} (\gamma'(y))^2 \,\mathrm{d}y.
\end{eqnarray*}
Since the second term above is finite by (C3), Lemma \ref{lemma4}
will follow if we can show that the first term goes to 0 as $n
\rightarrow
\infty$. Note that
\begin{eqnarray*}
&&n \int_0^1 (g_n(y))^2 \bigl(y(1-y)\bigr)^{-3/2} \,\mathrm{d}y \\
&&\quad\le 2^{3/2} n \biggl(\int_0^{{1}/{n}} \sqrt{y}   \,\mathrm{d}y +
\int_{1-{1}/{n}}^1 \sqrt{1-y}   \,\mathrm{d}y \biggr)
+ \frac{1}{n} \int_{{1}/{n}}^{1-1/n} y^{-3/2} (1-y)^{-3/2} \,\mathrm{d}y \\
&&\quad\le \frac{8 \sqrt{2}}{3 \sqrt{n}} + (1-n^{-1})^{-3/4} n^{-1/4}
\int_{0}^{1} y^{-3/4} (1-y)^{-3/4} \,\mathrm{d}y \rightarrow0.
\end{eqnarray*}
\upqed\end{pf}

\begin{lemma}
\label{lemma5} Let $U_{n\dvtx i} $ denote the $i$th order statistic of an
i.i.d.~sample of size $n$ from the uniform distribution over $(0,
1)$. Define $\mathcal{A}_{m, n} = \bigcap_{1 \le i \le n} \{
U_{n\dvtx i}(1-U_{n\dvtx i}) > \mu_{n\dvtx  i} (1-\mu_{n\dvtx  i})/m \}$. We then have
$ {\lim_{m \rightarrow\infty}
\sup_{n \ge1} P(\mathcal{A}_{m, n}) =1}$.
\end{lemma}

\begin{pf} By symmetry considerations, we only need to prove
%
\begin{equation}
\lim_{m \rightarrow\infty} \sup_{n \ge1} P\biggl(\bigcap_{1 \le i \le n/2}
\{U_{n\dvtx i}(1- U_{n\dvtx i})
> \mu_{n\dvtx i} (1-\mu_{n\dvtx i}) /m\}\biggr)=1.
\label{eqll4}
\end{equation}

For any $\varepsilon>0$, we can choose $n_0$ such that for all $n>n_0$,
$P(U_{n\dvtx ((n+1)/2)}\ge2/3)<\varepsilon/2$ and
\begin{eqnarray*}
&&P\biggl(\bigcap_{1\le i \le n/2} \{U_{n\dvtx i}(1- U_{n\dvtx i})
> \mu_{n\dvtx i} (1-\mu_{n\dvtx i}) /m\}\biggr)\\
&&\quad\ge P\biggl(\bigcap_{1\le i \le n/2} \{U_{n\dvtx i}
> 3\mu_{n\dvtx i} /m\}\biggr)-P\bigl(\bigl\{U_{n\dvtx ((n+1)/2)}\ge2/3\bigr\}\bigr)\\
&&\quad\ge P\biggl(\bigcap_{1\le i \le n/2} \{U_{n\dvtx i}
> 3\mu_{n\dvtx i} /m\}\biggr)-\varepsilon/2.
\end{eqnarray*}
Obviously, we can find a constant $m_0$ such that for all $m>m_0$,
\[
\sup_{1\le n\le n_0}P\biggl(\bigcap_{1 \le i \le n/2} \{U_{n\dvtx i}(1-
U_{n\dvtx i})
> \mu_{n\dvtx i} (1-\mu_{n\dvtx i}) /m\}\biggr)>1-\varepsilon.
\]
If we can choose a constant $m_1$ such that for all $m>m_1$,
\[
\sup_{n>n_0}P\biggl(\bigcap_{1 \le i \le n/2} \{U_{n\dvtx i}
> 3\mu_{n\dvtx i} /m\}\biggr)\ge1-\varepsilon/2,
\]
then, for all $m>\max(m_0,m_1)$,
\[
\sup_{n\ge1}P\biggl(\bigcap_{1 \le i
\le n/2} \{U_{n\dvtx i}(1- U_{n\dvtx i}) > \mu_{n\dvtx i} (1-\mu_{n\dvtx i}) /m\}
\biggr)>1-\varepsilon.
\]
Therefore, the
proof of Lemma \ref{lemma5} reduces to establishing
that
%
\begin{equation}
\lim_{m\to\infty}\sup_{n>1}P\biggl(\bigcap_{1 \le i \le n/2} \{U_{n\dvtx i} > \mu
_{n\dvtx i} /m\}\biggr)=1.
\label{lemll41}
\end{equation}

Recall the representation formula for the order statistics from a
sequence of uniform random variables,
$U_{n\dvtx i} \stackrel{\mathrm{dist}}{=} S_i/S_{n+1}$, where
$e_1, \ldots, e_{n+1}$ are i.i.d.~exponentially distributed random
variables with $E(e_i)=1$ and $S_i = e_1 + \cdots+ e_i$. If
\[
M= \inf_{1 \le i \le n <\infty}   \frac{S_i/i}{S_{n+1}/(n+1)} >\frac1m,
\]
then, for all $1 \le i\le n/2$, we have
$\frac{S_i/i}{S_{n+1}/(n+1)}>1/m$. This, in turn, implies that, as
$m\to\infty$,
\[
\lim_{m \rightarrow\infty}   \sup_{n \ge1}
P\biggl(\bigcap_{1 \le i \le n/2}\biggl\{\frac{S_i/i}{S_{n+1}/(n+1)} >
\frac1m\biggr\} \biggr)
\ge\lim_{m \rightarrow\infty} P(M>1/m)= P(M>0).
\]
Since $S_n /n \mathop{\longrightarrow}^{\mathrm{a.s.}} 1$ as $n \rightarrow\infty$, we have
$P(M>0)=1$.
This implies (\ref{lemll41}) and hence Lemma~\ref{lemma5} follows.
\end{pf}

\begin{pf*}{Proof of Theorem \protect\ref{asymnormality}}
We write
\[
\frac{1}{\sqrt{n}} \sum_{i=1}^n \psi\bigl(U_{n\dvtx  i}^{(1)}, \ldots, U_{n\dvtx
i}^{(d)}\bigr)
- \sqrt{n} \bar{\gamma} = I_n + \epsilon_n = S_{n,1} +S_{n,2} +
\epsilon_n,
\]
where
\begin{eqnarray*}
I_{n} &= & n^{-1/2} \sum_{i=1}^n \bigl(\psi\bigl(U_{n\dvtx  i}^{(1)}, \ldots,
U_{n\dvtx  i}^{(d)}\bigr)
- \gamma(\mu_{n\dvtx  i}) \bigr), \\
S_{n, 1}& = & n^{-1/2} \sum_{j=1}^d \sum_{i=1}^n \bigl(U_{n\dvtx i}^{(j)}-\mu
_{n\dvtx  i}\bigr) \psi_j (\mu_{n\dvtx  i}), \\
S_{n, 2} &= & I_n - S_{n, 1}, \\
\epsilon_n &= & n^{-1/2} \sum_{i=1}^{n} \gamma(\mu_{n\dvtx  i}) - \sqrt{n}
\int_0^1 \gamma(x) \,\mathrm{d}x.
\end{eqnarray*}
By Lemma \ref{lemma4}, $\epsilon_n \rightarrow0$ as $n \rightarrow
\infty$. We shall now
show that $S_{n, 2} \to^{P} 0$ as $n \rightarrow\infty$.

Since $ \max\{|U_{n\dvtx i}^{(j)}-\mu_{n\dvtx  i}|\dvtx  1 \le i \le n, 1 \le j
\le d \} \mathop{\longrightarrow}^{\mathrm{a.s.}} 0$,
by (C4) we have
\[
|S_{n, 2}|I_{\mathcal{A}_{m, n}} \le \frac{K(m)}{\sqrt{n}}
\sum_{j,k=1}^d \sum_{i=1}^n \bigl|\bigl(U_{n\dvtx i}^{(j)}-\mu_{n\dvtx  i}\bigr)
\bigl(U_{n\dvtx i}^{(k)}-\mu_{n\dvtx  i}\bigr)\bigr| \bigl( 1 + |\tilde{\psi}_{j,k} (\mu_{n\dvtx  i})|\bigr).
\]
By condition (C3), Lemma \ref{lemma1} and the Cauchy--Schwarz inequality,
we obtain
\begin{eqnarray*}
&&\frac{1}{\sqrt{n}} \sum_{i=1}^n
E\bigl|\bigl(U_{n\dvtx i}^{(j)}-\mu_{n\dvtx  i}\bigr) \bigl(U_{n\dvtx i}^{(k)}-\mu_{n\dvtx  i}\bigr)\bigr|
\bigl( 1 + |\tilde{\psi}_{j,k} (\mu_{n\dvtx  i})|\bigr)\\
&&\quad\le\frac{1}{n^{3/2}} \sum_{i=1}^n \mu_{n\dvtx  i} (1-\mu_{n\dvtx
i})|\tilde{\psi}_{j,k} (\mu_{n\dvtx  i})|
+ \frac{1}{\sqrt{n}}:= J_1 + J_2 + J_3 + \frac{1}{\sqrt{n}},
\end{eqnarray*}
where $J_1=n^{-3/2} \sum_{ 1\le i \le\sqrt{n}}\mu_{n\dvtx  i}
(1-\mu_{n\dvtx  i})|\tilde{\psi}_{j,k} (\mu_{n\dvtx
i})|$ and $J_2$, $J_3$ are similarly defined over
$\sqrt{n} < i <n-\sqrt{n}$ and $ n- \sqrt{n} \le i \le n$,
respectively. We have
\begin{eqnarray*}
J_1 &\le& \frac{2}{n} \sum_{1 \le i \le\sqrt{n}}
\bigl(\mu_{n\dvtx  i} (1-\mu_{n\dvtx  i})\bigr)^{3/2} |\tilde{\psi}_{j,k} (\mu
_{n\dvtx  i})| \\
&\sim& 2 \int_{0}^{1/\sqrt{n}} \bigl(x(1-x)\bigr)^{3/2}
|\tilde{\psi}_{j,k} (x)| \,\mathrm{d}x \rightarrow0
\end{eqnarray*}
as $n \rightarrow\infty$. Similarly, $J_3 \rightarrow0$ as $n
\rightarrow\infty$. Also,
as $ n \rightarrow\infty$,
\[
J_2 \le \frac{1}{n^{5/4}} \sum_{\sqrt{n} < i < n-\sqrt{n}}
\bigl(\mu_{n\dvtx  i} (1-\mu_{n\dvtx  i})\bigr)^{3/2} |\tilde{\psi}_{j,k}
(\mu_{n\dvtx  i})| \rightarrow0.
\]
That is, we have shown that as $n \rightarrow\infty$, for any given large
$m$, $S_{n, 2} I_{\mathcal{A}_{m, n}} \to^{P} 0$. We can now
choose a sequence of $m=m_n \rightarrow\infty$ such that $S_{n, 2}
I_{\mathcal{A}_{m, n}} \to^{P} 0$ as $ n \rightarrow\infty$.
By Lemma \ref{lemma5}, $I_{\mathcal{A}_{m, n}^c}\to^{P} 0$
and hence $S_{n,2} I_{\mathcal{A}_{m, n}^c} \to^{P} 0$
as $m \rightarrow\infty$. Therefore, $S_{n,2} \to^{P} 0$.

Define $W_{j, \ell} (x)= I(U_{\ell}^{(j)} \le x) -x$ for $ 1 \le j \le
d$ and $1 \le
\ell\le n$. Observe that $W_{j, 1} $ is $W_{j}$ defined in the proof
of Lemma \ref{lemma2} and
that $\hat{F}_{n; j}^{-1} (\frac{i}{n}) = U_{n\dvtx i}^{(j)}$.
By Bahadur's representation of quantiles (see, e.g.,
\cite{BR} or \cite{Kiefer}),
\[
\sup_{ 0 < t < 1} |\hat{F}_{n; j}(t) -t + \hat{F}_{n; j}^{-1}(t)
-t |
=\mathrm{O}(n^{-3/4} \log n) \qquad  \mbox{a.s. } \mbox{for } 1 \le j \le d.
\]
Hence,
\[
S_{n,1} = \frac{1}{\sqrt{n}} \sum_{\ell=1}^n Z_{n, \ell} + \mathrm{o}(1) \qquad\mbox{a.s.},
\]
where, for each $n$,
\[
Z_{n, \ell} = \frac{1}{n} \sum_{i=1}^n \sum_{j=1}^d
W_{j, \ell}(i/n) \psi_j (\mu_{n\dvtx  i})
\]
are i.i.d.~random variables with mean zero and
\begin{eqnarray*}
\operatorname{Var}(Z_{n, 1})
&=& \sum_{j, k=1}^d \frac{1}{n^2} \sum_{h, i=1}^n
\operatorname{Cov}\bigl(W_j(h/n), W_k(i/n)\bigr) \psi_j(\mu_{n\dvtx  h}) \psi_k(\mu_{n\dvtx  i}) \\
&=& \sum_{j, k=1}^d \frac{1}{n^2} \sum_{h, i=1}^n
\bigl(G_{j, k}(h/n, i/n) -hi n^{-2} \bigr) \psi_j(\mu_{n\dvtx  h}) \psi_k(\mu_{n\dvtx
i}) \\
&\to& \sum_{j, k =1}^d \int_0^1\!\! \int_0^1 \bigl(G_{j,k}(x, y)-xy\bigr) \psi
_j(x)\psi_k(y)\,\mathrm{d}x\,\mathrm{d}y.
\end{eqnarray*}
Recall that $G_{j, k}$ is the joint distribution of $(U_1^{(j)}, U_1^{(k)})$
and that $G_{j, j}(x, y)=\min(x,y)$.
To establish the convergence above, fix $j, k$ and split the second sum
above into cases according to
whether $h$, or $i$, is: less than $\epsilon n$; between $\epsilon n$
and $(1-\epsilon)n;$ greater than
$(1-\epsilon)n$. For example, when we sum over $\epsilon n \le h, i \le
(1-\epsilon)n$, then
it converges to $\int_{\epsilon}^{1-\epsilon} \int_{\epsilon
}^{1-\epsilon} H(x,y) \,\mathrm{d}x \,\mathrm{d}y$,
where $H(x, y) = (G_{j,k}(x, y) - xy) \psi_j(x) \psi_k(y)$.
The sum over $1 \le h < \epsilon n$ and $ \epsilon n \le i \le
(1-\epsilon)n$ can be shown to converge
to $\int_{0}^{\epsilon} \int_{\epsilon}^{1-\epsilon} H(x,y) \,\mathrm{d}x \,\mathrm{d}y$,
which, from the method of proof of Lemma
\ref{lemma2} and condition (C3), can be shown to converge
to 0 as $\epsilon\rightarrow0$. Similar convergences hold for other
ranges of $h$ and $i$.

It is now easy to see that the limit above can be written in
the form of $\sigma^2$ as stated in Theorem~\ref{asymnormality}. Note
that $|Z_{n,1}| \le\sum_{j=1}^d \frac{1}{n} \sum_{i=1}^n
|\psi_j(\mu_{n\dvtx i})|$. If $(1/\sqrt n) \frac{1}{n} \sum_{i=1}^n
|\psi_j(\mu_{n\dvtx i})| \to0$ for $j=1,2,\ldots,d,$ then the Lindeberg--L\'evy
condition holds. To see this, note that
\[
\Biggl(\frac1n \sum_{i=1}^n |\psi_j(\mu_{n\dvtx i})| \Biggr)^2
\le\frac1n \sum_{i=1}^n \bigl(\mu_{n\dvtx i}(1-\mu_{n\dvtx i})\bigr)^{-3/2}
\frac1n \sum_{i=1}^n \bigl(\mu_{n\dvtx i}(1-\mu_{n\dvtx i})\bigr)^{3/2}(\psi_j(\mu
_{n\dvtx i}))^2.
\]

By (C3), it is enough to establish that
$I_n = \frac{1}{n^{2}} \sum_{i=1}^n (\frac{i}{n+1}(1 - \frac
{i}{n+1}))^{-3/2}
\to0$. Since
\[
I_n \le\frac{4(n+1)^{3/2}}{ n^{2}} \biggl(\sum_{1\le i \le(n+1)/2} i^{-3/2}
+ \sum_{(n+1)/2 \le i \le n} (n+1-i)^{-3/2}\biggr) \to0,
\]
we have, by the Lindeberg--L\'evy central limit theorem,
\[
\frac{1}{\sqrt{n}} \sum_{\ell=1}^n Z_{n, \ell} \mathop{\longrightarrow}^{\mathrm{dist}} N(0,
\sigma^2).
\]
Hence, $S_{n, 1} \mathop{\to}^{\mathrm{dist}} N(0, \sigma^2)$, which completes the
proof of Theorem \ref{asymnormality}.
\end{pf*}

\section{Examples and counterexamples}\label{sec4}

We give some examples to show
our results and counterexamples to illustrate that conditions~(C1) and
(C2) are necessary for Theorem \ref{asconvergence} to hold.

\begin{enumerate}[(1)]
\item[(1)] Let $Z$ be a random variable with a continuous distribution
function $F$. Let
$g_j, 1 \le j \le d$, be continuous monotonically increasing
functions. For each $1 \le j \le d$, suppose $X_1^{(j)}, X_2^{(j)},
\ldots$ are independent random variables having the same distribution as
$g_j(Z)$. Applying Theorem \ref{asconvergence} and assuming
necessary integrability conditions, we get, after changing the
variable $y=F(x)$,
\[
\frac{1}{n} \sum_{i=1}^n \phi\bigl(X_{n\dvtx i}^{(1)}, \ldots, X_{n\dvtx i}^{(d)}\bigr)
\mathop{\longrightarrow}^{\mathrm{a.s.}} E \phi(g_1(Z), \ldots, g_d(Z))\qquad  \mbox{as } n
\rightarrow\infty.
\]

\item[(2)]
Let $ (X_1^{(1)}, \ldots, X_1^{(d)})$, $(X_2^{(1)}, \ldots,
X_2^{(d)}), \ldots$ be independent random vectors having the same
distribution as $(U_1, \ldots, U_d)$, where the $U_j$'s are uniformly
distributed over $(0,1)$. Let $F_{j,k}$ be the joint distribution of
$U_j$ and $U_k$. Suppose $\phi\dvtx  (0, 1)^d \rightarrow\R$ is defined by
$\phi(x_1, \ldots, x_d) = x_1^{\alpha_1} x_2^{\alpha_2} \cdots
x_d^{\alpha_d}$, where $\alpha_j \ge1$. Let $M=\alpha_1 + \cdots
+\alpha_d$.
Then $\psi=\phi$, $\gamma(x)=x^M$ and $\psi_j(x) =\alpha_j x^{M-1}$ for
$1 \le j \le d$. We have
\[
\frac{1}{n} \sum_{i=1}^n \bigl(X_{n\dvtx i}^{(1)}\bigr)^{\alpha_1} \cdots
\bigl(X_{n\dvtx i}^{(d)}\bigr)^{\alpha_d}
\mathop{\longrightarrow}^{\mathrm{a.s.}} \frac{1}{M +1}
\]
and
\[
n^{-1/2} \Biggl( \sum_{i=1}^n \bigl(X_{n\dvtx i}^{(1)}\bigr)^{\alpha_1} \cdots
\bigl(X_{n\dvtx i}^{(d)}\bigr)^{\alpha_d}
- \frac{n}{M +1} \Biggr) \mathop{\longrightarrow}^{\mathrm{dist.}} N(0, \sigma^2),
\]
where
\[
\sigma^2 =\frac{2}{M^2} \sum_{1 \le j <k \le d} \alpha_j \alpha_k \operatorname{Cov}(U_j^M, U_k^M) + \frac{1}{(M+1)^2(2M+1)} \sum_{j=1}^d \alpha_j^2.
\]

\item[(3)] The study of the statistical properties when there is a
loss of association among paired data has attracted a lot of attention
in various contexts, such as the broken sample problem, file linkage
problem and record linkage. For example, DeGroot and Goel initiated
the investigation of estimating the correlation coefficient of a
bivariate normal
distribution based on a broken random sample in \cite{DG}. Copas and
Hilton proposed statistical models to measure the evidence that a
pair of records relates to the same individuals in \cite{CH}. Chan
and Loh considered an approximation of the likelihood computation
for large broken sample in \cite{CL}.
Bai and Hsing, in \cite{BH}, proved that there does not exist any consistent
discrimination rule for the correlation coefficient, $\rho$, between
$X$ and $Y$ when the paired sample is broken, that is, the association
between $X$ and $Y$ is lost.
When pairing is lost, the $X$'s and $Y$'s behave as if they were
independent as far as
first order asymptotics, such as the law of large numbers (see Theorem
\ref{asconvergence}), are concerned.
\end{enumerate}

\begin{example}  This example shows that condition (C1) is
necessary for Theorem \ref{asconvergence} to hold. Let
\[
\phi(x, y) =
\cases{
1, & \quad $\mbox{if } 0 < x =y <1,$ \cr
0, & \quad $\mbox{if } 0 < x \neq y <1.$
}
\]
Let $\{(X_i, Y_i)\dvtx  1 \le i \le n\} $ be a sequence of i.i.d.~random
vectors. We further suppose that $X_i$ and $Y_i$ are independent and
uniformly distributed over $(0,1)$. Since $\phi$ is bounded,
(C2) holds, whereas (C1) does not hold.
We further note that $P(X_{n\dvtx i} \neq Y_{n\dvtx i} ) =1$ for $1 \le i \le n$.
Hence, $\sum_{i=1}^n \phi(X_{n\dvtx i}, Y_{n\dvtx i}) =0$, but $\int_0^1
\phi(x, x) \,\mathrm{d}x=1$.
\end{example}

\begin{example} This example shows that condition (C2) is
necessary for Theorem \ref{asconvergence} to hold. Let $\tilde{S}_0
=(0, 1)^2$ and, for $m \ge1$, define $\tilde{S}_m =
(\frac{m}{m+1}, 1)^2$, $S_m = \tilde{S}_{m-1} \setminus
\tilde{S}_m.$ Let $L_m$ be the union of three line segments:
\begin{eqnarray*}
L_m &=& \biggl\{ \biggl(\frac{m}{m+1}, y \biggr)\dvtx  \frac{m}{m+1} \le y
\le1 \biggr\} \cup
\biggl\{ \biggl(x, \frac{m}{m+1} \biggr)\dvtx  \frac{m}{m+1} \le x \le1
\biggr\} \\
&&{}  \cup\biggl\{(x, x)\dvtx  \frac{m-1}{m} \le x \le\frac{m}{m+1}
\biggr\}.
\end{eqnarray*}
Let $C_m$ be the region inside $S_m$ which is distance $\epsilon_m$
within $L_m$, where $\epsilon_m$ is chosen so that the area of $C_m$ is
$m^{-8}$. Write $A_m = S_m \setminus C_m$. Let $\phi$ be a
continuous on $(0,1)^2$ satisfying $\phi=1$ on the
diagonal, $\phi= m^3$ on $A_m$ and $1 \le\phi\le m^3$ on $C_m$.

Let $\{U_i, V_j\dvtx  1 \le i,j \le n\}$ be independent and uniformly
distributed on $(0,1)$.
Define $W_n =(U_{n\dvtx n}, V_{n\dvtx n})$ and $a_n = \lceil n ^{1/2} \rceil$.
Observe that
\begin{eqnarray}\label{eq4.1}
\frac1n \sum_{i=1}^n \phi(U_{n\dvtx  i} , V_{i\dvtx n})
&\ge& \frac1n \phi(W_n)
\ge \frac1n \sum_{m \ge\sqrt{n}} m^3   I( W_n \in A_m )\nonumber\\
&\ge& \sqrt{n} \sum_{m \ge\sqrt{n}}\bigl( I( W_n \in S_m )
- I( W_{n} \in C_m )\bigr) \\
&\ge& \sqrt{n} \biggl( I( W_n \in\tilde{S}_{a_n} ) -
I\biggl( W_n \in \bigcup_{m \ge\sqrt{n}} C_m \biggr)\biggr).
\nonumber
\end{eqnarray}
We now claim that
%
\begin{eqnarray}
I( W_n \in\tilde{S}_{a_n} ) & \longrightarrow& 1   \qquad \mbox{a.s.}, \label{claim1} \\
I\biggl( W_n \in\bigcup_{m \ge\sqrt{n}} C_m \biggr) & \longrightarrow& 0
  \qquad \mbox{a.s.} \label{claim2}
\end{eqnarray}
as $n \rightarrow\infty$. To prove (\ref{claim1}), observe that
\[
P(W_n \notin\tilde{S}_{a_n} )\le2 P\bigl(U_{n\dvtx n} \le a_n/(1+a_n)\bigr) =2
\biggl(1 -\frac{1}{1+a_n} \biggr)^n \approx2 \mathrm{e}^{-\sqrt{n}}.
\]
This yields
\[
\sum_{n =1}^{\infty} P(W_n \notin\tilde{S}_{a_n} ) < \infty,
\]
which, by the Borel--Cantelli lemma, implies that
$I( W_n \notin\tilde{S}_{a_n} \mbox{ i.o.})=0$, proving (\ref
{claim1}). To prove~(\ref{claim2}), it suffices to show that
%
\begin{equation}
\label{BC1} \sum_{n =1}^{\infty} P\biggl( W_n \in\bigcup_{m \ge
\sqrt{n}} C_m \biggr) < \infty.
\end{equation}
We again consider the $n$th term in the series in (\ref{BC1}):\vspace*{-2pt}
\begin{eqnarray*}
P\biggl( W_n \in\bigcup_{m \ge\sqrt{n}} C_m \biggr)
&\le& \sum_{m \ge\sqrt{n}} P(W_n \in C_m )\\[-1pt]
&\le& \sum_{m \ge\sqrt{n}} P\bigl((U_i, V_j) \in C_m \mbox{ for some } 1
\le i, j \le n \bigr)\\[-1pt]
&\le& n^2 \sum_{m \ge\sqrt{n}} P\bigl((U_1, V_1) \in C_m \bigr)= n^2 \sum_{m \ge\sqrt{n}} m^{-8}
\le C n^{-3/2}
\end{eqnarray*}\vspace*{-1pt}
and hence the infinite series in (\ref{BC1}) is finite. This
completes the proof of (\ref{claim2}). Thus, by (\ref{eq4.1}), $\frac
{1}{n} \sum_{i=1}^n
\phi(U_{n\dvtx  i} , V_{i\dvtx n})$ diverges. Furthermore, it is easy to see
that condition (C2) does not hold. If (C2) were satisfied, that
would imply boundedness of $\gamma$ over $(1-c_0, 1)^2$, which
is not the case. This completes the construction of the
counterexample.
\end{example}\vspace*{-7pt}

\section*{Acknowledgements}\vspace*{-4pt} The authors wish to thank the Editor,
an Associate Editor
and the referees for useful suggestions which improved the paper. They
acknowledge helpful discussions with
Dr. Alex Cook. The research of G.J.~Babu was supported in
part by NSF Grants AST-0707833 and AST-0434234. The research of
Z.D.~Bai was supported in part by the National University of Singapore ARF
Grant R-155-000-079-112. The research of K.P.~Choi is supported in
part by the National University of Singapore ARF Grants
R-155-000-051-112 and R-155-000-102-112.

\printhistory


\vspace*{-5pt}


\begin{thebibliography}{99}
\bibitem{BR}
Babu, G.J. and Rao, C.R. (1988). Joint asymptotic distribution of
marginal quantiles and quantile functions in samples from a
multivariate population. \textit{J. Multivariate Anal.} \textbf{27} 15--23.
\MR{0971169}

\bibitem{BH}
Bai, Z.D. and Hsing, T. (2005).
The broken sample problem. \textit{Probab.
Theory Related Fields}
\textbf{131} 528--552.
\MR{2147220}

\bibitem{Bill} Billingsley, P. (1999). \textit{Convergence of Probability
Measures},
2nd ed. \textit{Wiley Series in Probability and Statistics: Probability and Statistics}.
New York: Wiley.
\MR{1700749}

\bibitem{CH} Copas, J.B. and Hilton, F.J. (1990). Record linkage: Statistical
models for matching computer records.
\textit{J. R. Statist. Soc. A} \textbf{153} 287--320.

\bibitem{CL}
Chan, H.P. and Loh, W.L. (2001). A file linkage problem of DeGroot
and Goel revisited.
\textit{Statist. Sinica} \textbf{11} 1031--1045.
\MR{1867330}

\bibitem{David} David, H.A. (1981). \textit{Order Statistics}. New York: Wiley.
\MR{0286226}


\bibitem{DG} DeGroot, M.H. and Goel, P.K. (1980). Estimation of the correlation
coefficient from a broken sample.
\textit{Ann. Statist.} \textbf{8} 264--278.
\MR{0560728}

\bibitem{HLP} Hardy, G.H., Littlewood, J.E. and P\'{o}lya, G. (1952).
\textit{Inequalities}. Cambridge: Cambridge Univ. Press.

\bibitem{Kiefer} Kiefer, J. (1970). Deviations between the sample quantile
process and the
sample d.f. In \textit{Nonparametric Techniques in Statistical Inference
(Proc. Sympos., Indiana Univ., Bloomington, Ind., 1969)} 299--319.
London: Cambridge Univ. Press.
\MR{0277071}

\bibitem{VM} Mangalam, V. (2010).
Regression under lost association. To appear.
\end{thebibliography}
\end{document}